\documentstyle[12pt]{amsart}

\newcommand{\bZ}{\Bbb Z}
\newcommand{\bR}{\Bbb R}
\newcommand{\bC}{\Bbb C}
\newcommand{\G}{\Gamma}
\newcommand{\blowup}{\overline{\bC P^2}}

\newtheorem{thm}{Theorem}

\newtheorem{lem}{Lemma}

\theoremstyle{remark}

\title{On the blowups of numerical Godeaux surfaces}
\author{Dieter KOTSCHICK}
\address{Mathematisches Institut, Universit\'e de B\^ale,
Rheinsprung 21, 4051 B\^ale, Suisse}

\begin{document}

\maketitle

\bigskip\bigskip

\centerline{({\it to appear in C.R. Acad. Sci. Paris})}

\bigskip\bigskip

\noindent
{\bf Abstract:} {\sl We give a short proof of the following result:
Let $X$ be a complex surface of general type. If the canonical divisor of
the minimal model of $X$ has selfintersection $= 1$, then $X$ is not
diffeomorphic to a rational surface.

Our proof is the natural extension of the argument given in \cite{Invent}
for the case when $X$ is minimal. This argument also gives information about
the
non--existence of certain smooth embeddings of $2$--spheres in $X$, if $X$
has geometric genus zero.}

\bigskip

\bigskip

\centerline{\bf Sur les \'eclat\'ees des surfaces num\'eriquement de Godeaux}

\bigskip
\noindent
{\bf R\'esum\'e:} {\sl On donne une d\'emonstration rapide du r\'esultat
suivant:
Soit $X$ une surface complexe de type g\'en\'eral. Si l'autointersection
du diviseur canonique du mod\`ele minimal de $X$ est $1$, alors $X$
n'est pas diff\'eomorphe \`a une surface rationnelle.

Notre d\'emonstration est l'extension naturelle de la m\'ethode utilis\'ee
dans \cite{Invent} au cas o\`u $X$ est minimal. Cet argument donne aussi
des informations sur la non-existence de certains plongements lisses de
sph\`eres dans $X$, si le genre g\'eom\'etrique de $X$ est nul.}

\bigskip\bigskip

\noindent
{\bf Version fran\c{c}aise abr\'eg\'ee:}
On connait un seul type de d\'eformation de surfaces complexes minimales de
type
g\'en\'eral qui sont hom\'eomorphes \`a des surfaces rationnelles, ce sont les
surfaces de Barlow \cite{Barlow}. Elles sont simplement connexes mais ont
la m\^eme cohomologie rationnelle que la surface de Godeaux, obtenue comme
quotient d'une surface de degr\'e $5$ dans $\bC P^3$ par un groupe d'ordre $5$.
On dit que de telles surfaces sont {\sl num\'eriquement de Godeaux}.

Dans \cite{Invent}, on a associ\'e \`a certaines vari\'et\'es lisses
un entier $\phi$, invariant par diff\'eo-\linebreak morphisme.
On y a montr\'e $\vert\phi\vert\geq 4$ pour les surfaces de Barlow
et $\phi =0$ pour les surfaces rationnelles qui sont hom\'eomorphes \`a une
surface de Barlow.
Par cons\'equent, une surface de Barlow et une surface rationelle ne sont pas
diff\'eomorphes.
Pidstrigach \cite{Pid}
a \'etendu le calcul de \cite{Invent}, et a montr\'e
$\vert\phi\vert\geq 2$ pour les surfaces simplement connexes et
num\'eriquement de Godeaux.

Dans cette Note, nous g\'en\'eralisons ces arguments aux \'eclat\'ees.

\medskip
\noindent
{\bf Th\'eor\`eme 1.} {\it Aucune \'eclat\'ee d'une surface num\'eriquement
de Godeaux n'est diff\'eo-\linebreak morphe \`a une surface rationnelle.}

\medskip

Signalons qu'il existe deux autres d\'emonstrations de ce th\'eor\`eme
\cite{pos}, \cite{FQ}, qui sont beaucoup plus compliqu\'ees que la
d\'emonstration
de cette Note. La d\'emonstration de \cite{FQ} entra\^ine un resultat plus
g\'en\'eral: aucune surface de type g\'en\'eral n'est diff\'eomorphe
\`a une surface rationnelle. D'autre part, la d\'emonstration que l'on donne
ici pr\'esente l'avantage d'entrainer
aussi le th\'eor\`eme suivant, plus fort que le r\'esultat de \cite{FQ}:

\medskip
\noindent
{\bf Th\'eor\`eme 2.} {\it Soient $Y$ une surface num\'eriquement de Godeaux et
$F\in H^2(Y,\bZ)$ une classe d'autointersection $-1$. Si $X$ est une
\'eclat\'ee de $Y$, alors il est impossible de repr\'esenter
$F\in H^2(Y,\bZ)\subset H^2(X,\bZ)$ par une sph\`ere diff\'erentiablement
plong\'ee dans $X$}.

\medskip
\noindent
Ce dernier r\'esultat \'etait d\'ej\`a connu dans le cas $X=Y$,
(corollaire 7.9 de \cite{PLMS}). On sait aussi que $\phi(Y)\neq 0$
entraine que $Y$ ne contient pas de sph\`eres diff\'erentiablement
plong\'ees d'autointersection nulle mais non-triviales en homologie
rationnelle,
d'apr\`es le th\'eor\`eme 4.2 de \cite{KM} et la remarque qui le suit.

\bigskip
\bigskip\bigskip

The only known minimal complex surfaces of general type homeomorphic to
rational
surfaces are the Barlow surfaces \cite{Barlow}. They all have the same
deformation type as complex manifolds and thus are diffeomorphic to each other.
The Barlow surfaces $B$ are simply connected
but have the same rational homology as the classical Godeaux surface obtained
as the quotient of a quintic in $\bC P^3$ by a group of order $5$. We call such
surfaces {\sl numerical Godeaux surfaces}.

In \cite{Invent} we introduced a numerical diffeomorphism invariant $\phi$
of certain $4$-- manifolds, and proved $\vert\phi (B)\vert\geq 4$
for the Barlow surfaces and $\phi (R)=0$ for the rational surfaces $R$
homeomorphic to them. Thus, the two
are not diffeomorphic. More precisely, we calculated $\vert\phi (B)\vert=8$ in
\cite{PLMS}, \cite{AmJour}. Pidstrigach \cite{Pid} extended the calculation
of \cite{Invent} to show that $\vert\phi (Y)\vert \geq 2$ for any simply
connected numerical Godeaux surface $Y$.

We now generalize this argument to blowups.

\begin{thm}
No blowup of a numerical Godeaux surface is diffeomorphic to a rational
surface.
\end{thm}
\begin{pf}
Let $Y$ be a simply connected numerical Godeaux surface, and $X$ the
$k$--fold blowup of $Y$. The differentiable structure of $X$ is uniquely
determined by that of $Y$.

Let $P\rightarrow X$ be the principal $SO(3)$--bundle with $w_2(P)=w_2(TX)$
and $p_1(P)=-3-k$. Associated with $P$ there is a generalized Donaldson
invariant
$\phi_{k,c}^X$. It is a map from a certain set of chambers in the positive
cone for the cup--product form on $H^2 (X,\bR)$ to the algebra of symmetric
polynomial functions of degree
$k$ on $H_2(X,\bZ)$. The class $c\in H^2(X,\bZ)$ is a lift of
$w_2(P)\in H^2(X,\bZ_2)$ and
$\phi$ depends on this choice only up to sign. See \cite{PLMS}, \cite{II} for
definitions and further properties.
We will work up to signs and denote the invariant by $\phi_k^X$.
(The numerical invariant $\phi$ used in the minimal case is the same as
$\phi_0$. If $k\leq 1$, there is only one chamber, cf. \cite{Invent},
\cite{PLMS}.)

Suppose $R$ is a rational surface homeomorphic to $X$. Then $R$ contains
$k+8$ disjointly smoothly embedded $2$--spheres of selfintersection $-1$. It
follows that their
cohomology classes $E_1, ... , E_{k+8}$ are linearly independent over $\bR$.
In the positive cone $H_+^2(R,\bR)$ there is a unique chamber $C_0$
containing the ray of classes orthogonal to all the $E_i$ in its interior.
This chamber is preserved by the self--diffeomorphisms of $R$ realizing
the reflections in the $E_i$ on (co)homology. By the ``orientation argument''
of \cite{PLMS}, $\phi_k^R(C_0)$ is divisible by all the $E_i$, see
Proposition (6.13) of \cite{PLMS}. Thus, $\phi_k^R(C_0)$ is a polynomial
of degree $k$ with $k+8$ linear factors. This means it is identically zero,
cf. Corollary (6.14) of \cite{PLMS}.

In the cohomology of $X$, let $F_1, ... , F_k$ be the classes of the
exceptional curves. Choose $8$ classes $F_{k+1}, ... , F_{k+8}$
of square $-1$ in $H^2(Y,\bZ) \subset H^2(X,\bZ)$, linearly independent
over $\bR$.
Then there is a unique chamber $C_1$ in the positive cone $H^2_+(X,\bR)$
containing the ray of classes perpendicular to all the $F_i$ in its
interior. As in the case of $R$ above, the polynomial $\phi_k^X(C_1)$
is divisible by the classes $F_i$, for $i\leq k$. (This is not true for
$i > k$, as shown below. Compare Theorem 2.)

Thus, $\phi_k^X(C_1)$ is a scalar multiple of the product $F_1 ... F_k$.
Through a repeated application of Theorem (6.9) of \cite{PLMS} and its
proof, we see that the scalar is a universal non--zero combinatorial factor
times $\phi_0^Y =\phi (Y)$. In particular, it follows from
$\phi (Y)\neq 0$ that $\phi_k^X(C_1)\neq 0$.

Combined with $\phi_k^R(C_0) = 0$, this implies that $X$ and $R$ are not
diffeomorphic because of the following Lemma.
\end{pf}

\begin{lem}
If $X$ and $R$ are diffeomorphic, then they are diffeomorphic by a
diffeomorphism $f\colon R \rightarrow X$ with $f^*(C_1)=C_0$.
\end{lem}

\noindent
This is a special case of a Theorem of Friedman and Qin \cite{FQ}. They
observed that such results can be obtained by elementary arguments, using
the analysis of Donaldson's $\G$--invariant \cite{Gamma} in \cite{FM}.
Only the {\sl existence} of the invariant is used. No
non--trivial {\sl calculation} is needed.

I first outlined the above proof of Theorem 1 in a letter to M. Kreck in
July 1988, cf. \cite{Ober}. At that time, it was incomplete because Lemma 1 was
not known.
Because of this, I gave a different proof \cite{pos} using the
$\G$--invariant \cite{Gamma}. Recently, Friedman and Qin \cite{FQ} announced a
proof of Theorem 1 using the invariants introduced in \cite{PLMS}, \cite{II}.
Unlike
the above proof, theirs uses polynomial invariants of high degree even
for the minimal case. Both \cite{pos} and \cite{FQ} depend on subtle
calculations with stable bundles on the minimal model,
and use an $SU(2)$ blowup formula to prove the result
for the blowups. The argument given here
uses only the by now routine stable bunde calculation carried out in
\cite{Invent}, \cite{Pid}, and the $SO(3)$ blowup formula obtained in
\cite{PLMS}.

The argument of Friedman and Qin \cite{FQ} has the advantage that it
extends to surfaces of general type for which the canonical class
of the minimal model has selfintersection $>1$, thus showing that
no surface of general type is diffeomorphic to a rational one.
It would be interesting to know whether there exist examples of surfaces
of general type to which this applies, beyond the examples covered
by Theorem 1.

The proof of Theorem 1 given here furnishes more information about embedded
spheres in blowups of numerical Godeaux surfaces than either \cite{pos}
or \cite{FQ}. Indeed, the author and G. Mati\'c showed in \cite{KM},
Theorem 4.2 and the subsequent
Remark, that the non--vanishing of $\phi (Y)$ for a numerical Godeaux surface
$Y$ implies that it contains no smoothly embedded spheres of selfintersection
zero representing nontrivial rational homology classes. It also implies
that $Y$ contains no smoothly embedded sphere of selfintersection $-1$,
see Corollary (7.9) in \cite{PLMS}. This generalizes as follows:

\begin{thm}
Let $Y$ be a numerical Godeaux surface and $F\in H^2(Y,\bZ)$ a class of
square $-1$.
If $X$ is a blowup of $Y$, then $F\in H^2(Y,\bZ)\subset
H^2(X,\bZ)$ can not be represented by a smoothly embedded sphere in $X$.
\end{thm}
\begin{pf}
In the proof of Theorem 1 we can take $F$ as one of the $F_i$ with
$i>k$. If $F$ is represented by a smoothly embedded sphere, then
the argument there shows that $\phi_k^X(C_1)$ has $k+1$ linear
factors and so is identically zero. This is a contradiction.
\end{pf}
Note that Lemma 1 is not needed in the proof of Theorem 2.
The argument in \cite{FQ} only proves Theorem 2 under the additional hypothesis
$F.K_Y = 1$.

\bigskip

\bibliographystyle{amsplain}

\end{document}